\documentclass[a4paper,12pt]{article}

\usepackage{amsmath,amssymb}
\allowdisplaybreaks

\begin{document}

\begin{center}

\textbf{\Large Commutative hypergroups  associated} 

\medskip

\textbf{\Large  with a hyperfield}

\medskip

\vspace{0,2 cm}

{\large Herbert Heyer, Satoshi Kawakami, Tatsuya Tsurii }

\medskip

{\large and Satoe Yamanaka}

\medskip


\end{center}

\begin{abstract}

Let $H$ be a  commutative hypergroup  and 
$L$ a discrete commutative hypergroup. 
In the present paper we introduce a commutative hypergroup 
$\mathcal{K}(H, \varphi, L)$ associated with a hyperfield $\varphi$ of $H$ 
based on $L$. Moreover for the hyperfield $\varphi$ of a compact 
commutative hypergroup $H$ of strong type based on a discrete commutative hypergroup 
$L$ of strong type, we introduce the dual hyperfield 
$\hat{\varphi}$ of $\hat{L}$ based on $\hat{H}$ and show that 
$\widehat{\mathcal{K}}(H, \varphi, L) \cong \mathcal{K}(\hat{L}, \hat{\varphi}, \hat{H})$.

\bigskip

\medskip



\end{abstract}

\bigskip

\textbf{\large 1. Introduction}

\bigskip

In a previous paper ([HKTY2]) we discussed the hypergroup structure 
of the space $\mathcal{K}(\hat{H} \cup \widehat{H_0}, \mathbb{Z}_q(2))$ 
under the assumptions that $H$ is a compact hypergroup of strong type and $H_0$ 
is a closed subhypergroup of $H$ with finite index $[H : H_0]$, 
$\hat{H}$ and $\widehat{H_0}$ denoting the dual of $H$ and $H_0$ respectively, 
and $\mathbb{Z}_q(2)$ signifying the $q$-deformation of $\mathbb{Z}_2$ 
with $0 < q \leq 1$ in the sense of [KTY]. 
It remained an open problem to investigate the hypergroup structure 
of the dual $\hat{\mathcal{K}}(\hat{H} \cup \widehat{H_0}, \mathbb{Z}_q(2))$  
of $\mathcal{K}(\hat{H} \cup \widehat{H_0}, \mathbb{Z}_q(2))$ 
at least under additional assumptions on $H$ and $H_0$. 
It  required a different approach to master that problem. 
The present solution to the problem relies on a generalization of the 
notion of a hyperfield, which had been successfully applied 
in the case of finite hypergroups in [HKKK]. 

\medskip

Let us emphasize that in the present discussion the dual 
hypergroup structure of $\mathcal{K}(\hat{H} \cup \widehat{H_0}, \mathbb{Z}_q(2))$ 
will be established for commutative hypergroups $H$ without 
assuming compactness of $H$ and finiteness of the index $[H:H_0]$ of 
$H_0$ in $H$. 

\medskip

The preliminary knowledge on hypergroups needed in the sequel 
can be taken from the traditional sources [BH] and [J]; 
some additional references on the structure of hypergroups to be 
consulted are [HK1], [HK2] and [HKY]. 

\medskip

In section 3 the appropriate generalization of the hyperfield 
method is presented. We are considering hyperfields $\varphi$ 
mapping elements $\ell$ of a countable discrete commutative hypergroup 
$L$ to compact subhypergroups $H(\ell)$ of a commutative hypergroup $H$. 
Then the space $\mathcal{K}(H, \varphi, L)$ is introduced and 
shown to be a commutative hypergroup (Theorem 3.1). 
Next we define the hypergroup $\mathcal{K}(\hat{L}, \hat{\varphi}, \hat{H})$ 
by the dual field of $\hat{\varphi}$ of $\varphi$ and 
obtain that for a compact commutative hypergroup $H$ of strong type   
the desired dual $\hat{\mathcal{K}}(H, \varphi, L)$ 
 is isomorphic to $\mathcal{K}(\hat{L}, \hat{\varphi}, \hat{H})$ 
 (Theorem 3.5). 
Natural conditions yield the identification of 
$\mathcal{K}(\hat{L}, \hat{\varphi}, \hat{H})$ with $\hat{L} \vee \hat{H}$ 
and with the substitution $S(Q \times L : Q \longrightarrow H)$ 
introduced by Voit in [V1]. 

\medskip

Section 4 contains various examples of the hyperfield method, 
in particular the extension of Voit's result in [V2] to 
higher dimensional tori.

\medskip

In section 5 we drop the assumption of compactness of $H$ 
and identify $\mathcal{K}(\hat{H} \cup \widehat{H_0}, \mathbb{Z}_q(2))$ 
with $\mathcal{K}(\hat{H}, \varphi, \mathbb{Z}_q(2))$ (Theorem 5.2) 
using the character theory for induced representations developed in [HKY].

\bigskip

\textbf{\large 2. Preliminaries} 

\medskip

For a locally compact space $X$ we shall mainly consider the subspaces $C_c(X)$ and $C_0(X)$ of the space $C(X)$ 
of  continuous functions on $X$ which have compact support or vanish at infinity respectively. 
By $M(X)$, $M^b(X)$ and $M_c(X)$ we abbreviate the spaces of all (Radon) measures on $X$, the 
bounded measures and the measures with compact support on $X$ respectively. 
Let $M^1(X)$ denote the set of probability measures on $X$ and $M^1_c(X)$ its subset $M^1(X) \cap M_c(X)$.
The symbol $\delta_x$ stands for the Dirac measures in $x \in X$. 

\medskip

A {\it hypergroup} $(K,*)$ is a locally compact space $K$ together with a {\it convolution} $*$ in $M^b(K)$ such 
that $(M^b(K),*)$ becomes a Banach algebra and that the following properties are fulfilled. 

\medskip

(H1) The mapping 
$$
(\mu, \nu )\longmapsto \mu *  \nu
$$

\ \ \ \ \ \ \ from $M^b(K) \times M^b(K)$ into $M^b(K)$ is continuous with respect to the 

\ \ \ \ \ \ \ weak  topology in $M^b(K)$. 

\medskip

(H2) For $x,y \in K$ the convolution $\delta_x * \delta_y$ belongs to $M^1_c(K)$.

\medskip

(H3) There exist a {\it unit} element $e \in K$ with 
$$
\delta_e * \delta_x = \delta_x * \delta_e = \delta_x 
$$

\ \ \ \ \ \ \ for all $x \in K$, and an {\it involution} 
$$
x \longmapsto x^-
$$
   
\ \ \ \ \ \ \ in $K$ such that 
$$
\delta_{x^-} * \delta_{y^-} = (\delta_y * \delta_x)^-
$$

\ \ \ \ \ \ \ and 
$$
e \in \text{supp}(\delta_x * \delta_y)~~\text{if and only if}~~x = y^-
$$

\ \ \ \ \ \ \ whenever $x,y \in K$. 

\medskip

(H4) The mapping 
$$
(x,y) \longmapsto \text{supp}(\delta_x * \delta_y)
$$

\ \ \ \ \ \ \ from $K \times K$ into the space $\mathcal{C}(K)$ of all compact subsets of $K$ furnished  

\ \ \ \ \ \ \ with Michael topology is continuous. 

\medskip

A hypergroup $(K, *)$ is said to be {\it commutative} if the convolution $*$ is commutative. 
In this case $(M^b(K), *, -)$ is a commutative Banach $*$-algebra with identity $\delta_e$. 
There is an abundance of hypergroups and there are various constructions (polynomial, Sturm-Liouville) as the reader 
may learn from the pioneering papers on the subject. 
\medskip

Let $(K,*)$ and $(L, \circ)$ be two hypergroups with units $e_K$ and $e_L$ respectively. A continuous mapping 
$\varphi:K \rightarrow L$ is called a hypergroup {\it homomorphism} if $\varphi(e_K)=  e_L$ and $\varphi$ is the unique linear, weakly continuous 
extension from $M^b(K)$ to $M^b(L)$ 
such that 
$$
\varphi(\delta_x) = \delta_{\varphi(x)}, \ 
\varphi(\delta_x^-) = \varphi(\delta_x)^- \  \textrm{and} \  
\varphi (\delta_x * \delta_y) = \varphi (\delta_x) \circ \varphi (\delta_y)
$$
whenever $x,y \in K$. 
If $\varphi:K \rightarrow L$ is also a homeomorphism, it will be called 
an {\it isomorphism} from $K$ onto $L$. 
An isomorphism from $K$ onto $K$ is called an {\it automorphism} of $K$. 
We denote by Aut($K$) the set of all automorphisms of $K$. 
Then Aut($K$) becomes a topological group equipped with the weak topology 
of $M^b(K)$. 
We call $\alpha$ an \textit{action} of a locally compact group $G$ on a hypergroup 
$H$ if $\alpha$ is a continuous homomorphism from $G$ into Aut($H$). 
Associated with the action $\alpha$ of $G$ on $H$ one can define a 
semi-direct product hypergroup $K=H\rtimes_{\alpha}G$.  

\medskip

\medskip

If the given hypergroup $K$ is commutative, its dual 
$\widehat{K}$ can be introduced as the set of all bounded continuous functions 
$\chi \not = 0$ on $M^b(K)$ satisfying
\begin{align*}
&\chi(\delta_x * \delta_y) = \chi(\delta_x) \chi(\delta_y) ~~{\rm and}~~ 
\chi(\delta_x^-) = \overline{\chi(\delta_x)}
\end{align*}
for all $x,y \in K$. This set of characters $\widehat{K}$ of $K$ 
becomes a locally compact space with  respect to 
the topology of uniform convergence on compact sets, 
but generally fails to be a hypergroup. 
If $\widehat{K}$ is 
a hypergroup, then $K$ is called a strong hypergroup or 
a hypergroup of strong type. If the dual 
$\widehat{K}$ of a strong hypergroup $K$ is also strong and 
$\widehat{\widehat{K}} \cong K$ holds, then $K$ is called a {\it Pontryagin} hypergroup 
or a hypergroup of  {\it Pontryagin type}.

\bigskip

\textbf{\large 3. Hyperfields and hypergroups}

\medskip

Let $H = (H, M^b(H), *,{}^-)$ be a commutative hypergroup with unit $h_0$ and 
$L = (L, M^b(L), \bullet, {}^-) = \{\ell_0, \ell_1, \cdots, \ell_n, \cdots\}$ 
a countable discrete commutative hypergroup 
 where $\ell_0$ is unit of $L$. 

\medskip
\noindent
{\bf Definition } For each $\ell \in L$, let $H(\ell)$ be a compact subhypergroup of 
$H$ satisfying the following conditions. 

\medskip

(1) $H(\ell_0) = \{h_0\}$ and $H(\ell^-) = H(\ell)$. 

\medskip

(2) $[H(\ell_i) * H(\ell_j)] \supset H(\ell_k)$ for 
$\ell_k \in {\rm supp}(\varepsilon_{\ell_i} \bullet \varepsilon_{\ell_j})$, 
where $[H(\ell_i) * H(\ell_j)]$ is the  

~~~~~compact subhypergroup of $H$ generated by $H(\ell_i)$ 
and $H(\ell_j)$. 

\medskip
\noindent
Then we call 
$$
\varphi : L \ni \ell \longmapsto H(\ell) \subset H
$$
a {\it hyperfield} of $H$ based on $L$. 

\medskip

We denote the normalized Haar measure $\omega_{H(\ell)}$ of $H(\ell)$ by $e(\ell)$ and  
 note that  condition (2) implies 

\medskip

(3) $e(\ell_i) * e(\ell_j) * e(\ell_k) = e(\ell_i) * e(\ell_j)$ for 
$\ell_k \in {\rm supp}(\varepsilon_{\ell_i} \bullet \varepsilon_{\ell_j})$. 

\medskip

Putting 
$$\mathcal{K}(H, \varphi, L) := \{ (\delta_h*e(\ell)) \otimes \varepsilon_{\ell} 
\in M^b(H) \otimes M^b(L) : h \in H, \ell \in L \}$$
 
\medskip
\noindent
one sees  that 
$$
\mathcal{K}(H, \varphi, L) = Q(\ell_0) \cup Q(\ell_1) \cup \cdots \cup Q(\ell_n) \cup \cdots
$$
is a locally compact  space, where 
\begin{align*}
&Q(\ell) := H/H(\ell)~~{\rm for~all~}\ell \in L,\\
&Q(\ell_0) := H. 
\end{align*}

Now we shall introduce a convolution $\circ$ and an involution ${}^-$ on 
$\mathcal{K}(H, \varphi, L)$ as elements of $M^b(H) \otimes M^b(L)$ in 
order to obtain the following theorem generalizing a result in  [HKKK]. 

\medskip
\noindent
{\bf Theorem 3.1 } Let $\varphi$ be a hyperfield of a 
commutative hypergroup $H$ based on a countable discrete 
commutative hypergroup $L$. Then $\mathcal{K}(H, \varphi, L)$ is a commutative hypergroup. 

\medskip
\noindent
{\bf Proof } The set $\{\delta_h * e(\ell) : h \in H\}$ is a commutative 
hypergroup isomorphic to the quotient hypergroup 
$$
Q(\ell) = H/H(\ell)~~{\rm for~all~}\ell \in L.
$$ 
We see that 
$$
M^b(\mathcal{K}(H, \varphi, L)) =  \sum_{\ell \in L} \!{}^{\oplus} \ \! M^b(Q(\ell)). 
$$
Next we examine the convolution on $\mathcal{K}(H, \varphi, L)$ in detail. 
For $\ell_i, \ell_j \in L$ we denote 
the set 
$ \{k : \ell_k \in {\rm supp}(\varepsilon_{\ell_i} \bullet \varepsilon_{\ell_j}) \}$ 
by $s(\ell_i, \ell_j)$. Given $h_p, h_q \in H$ and $\ell_i, \ell_j \in L$,
\begin{align*}
&((\delta_{h_p} * e(\ell_i)) \otimes \varepsilon_{\ell_i}) \circ 
((\delta_{h_q} * e(\ell_j)) \otimes \varepsilon_{\ell_j})\\
=& (\delta_{h_p} * \delta_{h_q} * e(\ell_i) * e(\ell_j)) \otimes 
(\varepsilon_{\ell_i} \bullet \varepsilon_{\ell_j})\\
=& (\delta_{h_p} * \delta_{h_q} * e(\ell_i) * e(\ell_j) * e(\ell_k)) \otimes 
\left( \sum_{k \in s(\ell_i, \ell_j)} n_{ij}^k \varepsilon_{\ell_k} \right)\\
=& \sum_{k \in s(\ell_i, \ell_j)} n_{ij}^k 
(\delta_{h_p} * \delta_{h_q} * e(\ell_i) * e(\ell_j) * e(\ell_k)) \otimes 
\varepsilon_{\ell_k}
\end{align*}
by condition (3) derived from the defining properties of the hyperfield $\varphi$. 
In conclusion the convolution in $\mathcal{K}(H, \varphi, L)$ is well-defined, and 
its associativity holds. 

\medskip
Now we note that 
$$
{\rm supp}(((\delta_{h_p} * e(\ell_i)) \otimes \varepsilon_{\ell_i}) \circ 
((\delta_{h_q} * e(\ell_j)) \otimes \varepsilon_{\ell_j}))) = \bigcup_{k \in s(\ell_i, \ell_j)} 
((h_p * h_q *[H(\ell_i) * H(\ell_j)])/H(\ell_k))
$$
is compact. 

\medskip

In order to verify the defining property of the involution ${}^-$ 
of $\mathcal{K}(H, \varphi, L)$ we compute 
\begin{align*}
&((\delta_{h_p} * e(\ell_i)) \otimes \varepsilon_{\ell_i}) \circ 
((\delta_{h_q} * e(\ell_j)) \otimes \varepsilon_{\ell_j})^-\\
=& ((\delta_{h_p} * e(\ell_i)) \otimes \varepsilon_{\ell_i}) \circ 
((\delta_{h_q}^- * e(\ell_j)) \otimes \varepsilon_{\ell_j}^-)\\
=& \sum_{k \in s(\ell_i, \ell_j^-)} n_{ij}^k 
(\delta_{h_p} * \delta_{h_q}^- * e(\ell_i) * e(\ell_j) * e(\ell_k)) \otimes 
\varepsilon_{\ell_k}.
\end{align*}
These equalities imply that 
$$
((\delta_{h_p} * e(\ell_i)) \otimes \varepsilon_{\ell_i})
= ((\delta_{h_q} * e(\ell_j)) \otimes \varepsilon_{\ell_j})^-
$$
holds if and only if 
$$
(h_0, \ell_0) \in {\rm supp}(((\delta_{h_p} * e(\ell_i)) \otimes \varepsilon_{\ell_i}) \circ 
((\delta_{h_q} * e(\ell_j)) \otimes \varepsilon_{\ell_j})). 
$$
The remaining axioms of a hypergroup are easily verified. 
All together  the desired conclusions are established.  
\hspace{\fill}[Q.E.D.]

\medskip
\noindent
{\bf Remark } If $\varphi(\ell) = H(\ell) = \{h_0\}$ for all $\ell \in L$, 
then $\mathcal{K}(H, \varphi, L)$ is the direct product hypergroup $H \times L$. 

\medskip

Now let $H$ be a compact commutative hypergroup of strong type such that 
$\hat{H}$ is a discrete commutative hypergroup and $L$ be a 
discrete commutative hypergroup of strong type such that $\hat{L}$ is a 
compact hypergroup. 
For $\chi \in \hat{H}$ we put 
$$
Y(\chi) := \{\ell \in L : \chi \in H(\ell)^\bot \}
$$
where 
$$
H(\ell)^\bot := \{\chi \in \hat{H} : \chi(h) = 1~~{\rm for~all}~~h \in H(\ell)\}
$$
denotes the annihilator of $H(\ell)$ in $\hat{H}$.

\medskip
\noindent
{\bf Lemma 3.2 } Given $\chi \in \hat{H}$, $Y(\chi)$ is a subhypergroup 
of $L$ such that $Y(\chi^-) = Y(\chi)$, and 
$Y(\chi_i) \cap Y(\chi_j) \subset Y(\chi_k)$ for $\chi_k \in {\rm supp}(\varepsilon_{\chi_i} * \varepsilon_{\chi_j})$. 

\medskip
\noindent
{\bf Proof } First of all we note that 
$$
H(\ell_i)^\bot \cap H(\ell_j)^\bot = [H(\ell_i) * H(\ell_j)]^\bot \subset H(\ell_k)^\bot
$$
for $\ell_k \in {\rm supp}(\varepsilon_{\ell_i} \bullet \varepsilon_{\ell_j})$ by 
the defining property (2) of the hyperfield $\varphi$ of $H$ based on $L$.

For $\ell_i, \ell_j \in Y(\chi)$, $\chi \in H(\ell_i)^\bot \cap H(\ell_j)^\bot \subset H(\ell_k)^\bot$ 
whenever $\ell_k \in {\rm supp}(\varepsilon_{\ell_i} \bullet \varepsilon_{\ell_j})$, i.e., 
$\chi \in H(\ell_k)^\bot$. Then 
${\rm supp}(\varepsilon_{\ell_i} \bullet \varepsilon_{\ell_j}) \subset Y(\chi)$.  
For $\ell \in Y(\chi)$ we have $\ell^- \in Y(\chi)$ by  property (1) of 
the definition of the hyperfield $\varphi$ which states $H(\ell^-) = H(\ell)$ for all $\ell \in L$. Then 
$Y(\chi)$ appears to be a subhypergroup of $L$. 
  
For $\ell \in Y(\chi_i) \cap Y(\chi_j)$ we see that 
$\chi_i \in H(\ell)^\bot$ and $\chi_j \in H(\ell)^\bot$.  
Since $H(\ell)^\bot$ is a subhypergroup of $\hat{H}$ we get 
$$
{\rm supp}(\varepsilon_{\chi_i} * \varepsilon_{\chi_j}) \subset H(\ell)^\bot, 
$$
which implies that $\chi_k \in H(\ell)^\bot$ for $\chi_k \in {\rm supp}(\varepsilon_{\chi_i} * \varepsilon_{\chi_j})$.
Altogether we arrive at the fact that $\ell \in Y(\chi_k)$, i.e., 
$Y(\chi_i) \cap Y(\chi_j) \subset Y(\chi_k)$. \hspace{\fill}[Q.E.D.]

\bigskip

Denoting $Y(\chi)^\bot$ by $\hat{L}(\chi)$ for $\chi \in \hat{H}$ we easily see  that 
$\hat{L}(\chi)$ is a closed subhypergroup of the compact hypergroup $\hat{L}$ 
satisfying properties (1) and (2) of the hyperfield $\varphi$ by Lemma 3.2. 
This leads to the following.

\medskip
\noindent
{\bf Lemma 3.3 }  The mapping 
$$
\hat{\varphi} : \hat{H} \ni \chi \longmapsto \hat{L}(\chi) \subset \hat{L}
$$
is a hyperfield of $\hat{L}$ based on $\hat{H}$.

\medskip
\noindent
{\bf Definition } The hyperfield $\hat{\varphi}$ is called the{ \it dual} of the 
hyperfield 
$$
\varphi : L \ni \ell \longmapsto H(\ell) \subset H
$$
of $H$ based on $L$. We note that the duality 
$\hat{\hat{\varphi}} = \varphi$ holds if $H$ and 
$L$ are Pontryagin. 

\bigskip

As a consequence of these preparations we obtain a commutative hypergroup 
$$
\mathcal{K}(\hat{L}, \hat{\varphi}, \hat{H}) = 
\{(\delta_{\rho} * e(\chi)) \otimes \varepsilon_{\chi} : \rho \in \hat{L}, \chi \in \hat{H} \},
$$
where $e(\chi)$ denotes the normalized Haar measure of $\hat{L}(\chi)$.

\medskip
 
The following statements are easily verified. 

\medskip
\noindent
{\bf Lemma 3.4 } 

\medskip

(i) For each $\chi \in \hat{H}$ and $\ell \in L$, $\ell \in Y(\chi)$ if and only if 
$\chi \in H(\ell)^\bot$

\medskip

(ii) For each $\chi \in \hat{H}$ and the Haar measure $e(\ell)$ of $H(\ell)$
\begin{eqnarray*}
\chi(e(\ell))=
\left\{ 
\begin{array}{ll}
1 &  ~~~{\rm if}~~\chi \in H(\ell)^\bot  \\
0 &  ~~~{\rm otherwise} \\
\end{array} 
\right.
\end{eqnarray*}

\medskip

(iii) For each $\ell \in L$ and the Haar measure $e(\chi)$ of $\hat{L}(\chi)$
\begin{eqnarray*}
e(\chi)(\ell)=
\left\{ 
\begin{array}{ll}
1 &  ~~~{\rm if}~~\ell \in Y(\chi)  \\
0 &  ~~~{\rm otherwise} \\
\end{array} 
\right.
\end{eqnarray*}

\medskip

(iv) For each  $\chi \in \hat{H}$ and $\ell \in L$, $\chi(e(\ell)) = e(\chi)(\ell)$.

\bigskip

\noindent
Now, we arrive at the dual version of the statement of Theorem 3.1. 

\medskip
\noindent
{\bf Theorem 3.5 } Let $\varphi$ be a hyperfield of a compact commutative 
hypergroup $H$ of strong type based on a discrete commutative hypergroup 
$L$ of strong type. Then 
$$
\hat{\mathcal{K}}(H, \varphi, L) \cong \mathcal{K}(\hat{L}, \hat{\varphi}, \hat{H}).  
$$
If $H$ and $L$ are Pontryagin hypergroups, then $\mathcal{K}(H, \varphi, L)$ is 
also Pontryagin. Moreover the sequence 
$$
1 \longrightarrow H \longrightarrow \mathcal{K}(H, \varphi, L) \longrightarrow 
L \longrightarrow 1
$$
is exact and the dual  sequence  

$$
1 \longrightarrow \hat{L} \longrightarrow \mathcal{K}(\hat{L}, \hat{\varphi}, \hat{H})
\longrightarrow  \hat{H} \longrightarrow 1
$$
is exact as well. In particular  $\mathcal{K}(H, \varphi, L)$ and $ \mathcal{K}(\hat{L}, \hat{\varphi}, \hat{H})$ 
are extension hypergroups of $L$ by $H$ and $ \hat{H}$ by  $\hat{L}$ respectively. 

\medskip
\noindent
{\bf Proof } Clearly 
$$
\hat{\mathcal{K}}(H, \varphi, L) \supset 
\mathcal{K}(\hat{L}, \hat{\varphi}, \hat{H}).
$$ 
It remains to be shown that 
$$
\hat{\mathcal{K}}(H, \varphi, L) \subset 
\mathcal{K}(\hat{L}, \hat{\varphi}, \hat{H}).
$$ 
Let $\tau$ be a character of $\mathcal{K}(H, \varphi, L)$. 
Then there exists $\chi \in \hat{H}$ such that 
$$
\tau((\delta_h * e(\ell)) \otimes \varepsilon_\ell) = \chi(h) \chi(e(\ell))\rho(\ell) 
= \rho(\ell)e(\chi)(\ell)\chi(h) = (\delta_{\rho} * e(\chi)) \otimes \varepsilon_{\chi})(\ell,h) 
$$
for some $\rho \in \hat{L}$ by Lemma 3.4. Consequently 
$$
\tau = (\delta_{\rho} * e(\chi)) \otimes \varepsilon_{\chi} \in \mathcal{K}(\hat{L}, \hat{\varphi}, \hat{H}).  
$$

The assertion concerning the Pontryagin property follows from the fact that 
the dual $\hat{\hat{\varphi}}$ of the dual hyperfield $\hat{\varphi}$ is $\varphi$. 

\medskip

Now let $e(H)$ denote the normalized Haar measure of the compact hypergroup $H$. Then 
\begin{align*}
Q  &:= \{(e(H) \otimes \varepsilon_{\ell_0}) \circ (\delta_h \otimes \varepsilon_{\ell_i}) : h \in H, \ell_i \in L\}\\
&= \{e(H) \otimes \varepsilon_{\ell_i} : \ell_i \in L \}
\end{align*}
is the quotient hypergroup $\mathcal{K}(H, \varphi, L)/H$ isomorphic to $L$. 
This means that $\mathcal{K}(H, \varphi, L)$ is an extension hypergroup 
of $L$ by $H$. 
\hspace{\fill}[Q.E.D.]

\bigskip
\medskip
\noindent
{\bf Remark } 
\medskip

(1) If $\varphi(\ell_0) = \{h_0\}$ and $\varphi(\ell) = H$ for all $\ell \in L$ such that 
$\ell \not = \ell_0$, then $\mathcal{K}(H, \varphi, L)$ is the 
hypergroup join $H \vee L$. 
Moreover, if $H$ and $L$ are strong, then 
$$
\mathcal{K}(\hat{L}, \hat{\varphi}, \hat{H}) = \hat{L} \vee \hat{H}. 
$$

(2) Let $Q := H/H_0$ for a closed subhypergroup $H_0$ of $H$. 
If $\varphi(\ell_0) = \{h_0\}$ and $\varphi(\ell) = H_0$ for all $\ell \in L$, 
$\ell \not = \ell_0$, then 
$$
\mathcal{K}(H, \varphi, L) = S(Q \times L : Q \longrightarrow H), 
$$
where the latter symbol denotes the substitution hypergroup obtained by 
substituting $Q$ in $Q \times L$ by $H$, in the sense 
of Voit [V1].  
We note that $H_0$ is not assumed to be  open which means that our definition 
is a generalization of Voit's substitution.

\bigskip

\textbf{\large 4. Examples of hyperfields}

\medskip

Let $\mathbb{Z}_q(2) = \{\ell_0, \ell_1\}$ be the hypergroup of order two, where the 
convolution structure is given by 
$$
\varepsilon_{\ell_1} \bullet \varepsilon_{\ell_1} 
= q \varepsilon_{\ell_0} + (1-q)\varepsilon_{\ell_1}
$$ 
for $0 < q \leq 1$.

\medskip
\noindent
{\bf Example 4.1 } If 
$H = \mathbb{T} = \{z \in \mathbb{C} : |z| = 1\}$, 
$L = \mathbb{Z}_q(2) = \{\ell_0, \ell_1\}$ and 
$\varphi : \mathbb{Z}_q(2) \ni \ell \longmapsto H(\ell) \subset H$ 
with $\varphi(\ell_0) = H(\ell_0) = \{1\}$, 
$\varphi(\ell_1) = H(\ell_1) = C_n := \{z \in \mathbb{T} : z^n = 1\}$ ($n \in \mathbb{N}$), 
then 
$$
\mathcal{K}(H, \varphi, L) = \mathcal{K}(\mathbb{T}, \varphi, \mathbb{Z}_q(2)) 
= \mathbb{T} \cup \mathbb{T}
$$
is a commutative hypergroup which coincides with  Voit's 
commutative hypergroup on the two tori 
$\mathbb{T} \cup \mathbb{T}$ ([V2]).

This means that $\mathcal{K}(\mathbb{T}, \varphi, \mathbb{Z}_q(2))$ determines the 
commutative hypergroup structure on $\mathbb{T} \cup \mathbb{T}$ with 
parameter $(n,q)$, $n \in \mathbb{N}$, $0 < q \leq 1$. 

Obviously $\mathcal{K}(\mathbb{T}, \varphi, \mathbb{Z}_q(2))$ is Pontryagin 
and 
$$
\hat{\mathcal{K}}(\mathbb{T}, \varphi, \mathbb{Z}_q(2)) = 
\mathcal{K}(\mathbb{Z}_q(2), \hat{\varphi}, \mathbb{Z}), 
$$
where the dual field $\hat{\varphi}$ of $\varphi$ is given as 
$$
\hat{\varphi} : \mathbb{Z} \ni k \longmapsto \hat{\varphi}(k) \subset \mathbb{Z}_q(2) 
$$
with 
\begin{eqnarray*}
\hat{\varphi}(k)=
\left\{ 
\begin{array}{ll}
\{\ell_0\} &  ~~~{\rm for}~~k \in n \mathbb{Z}  \\
\mathbb{Z}_q(2) &  ~~~{\rm otherwise}. \\
\end{array} 
\right.
\end{eqnarray*}

\bigskip

{\bf Example 4.2 } If  
$H = \mathbb{T} \times \mathbb{T} = \mathbb{T}^2$,  
$L = \mathbb{Z}_q(2) = \{\ell_0, \ell_1\}$ and 
$\varphi : \mathbb{Z}_q(2) \ni \ell \longmapsto H(\ell) \subset H$ 
with $\varphi(\ell_0) = H(\ell_0) = \{ (1,1) \}$, 
$\varphi(\ell_1) = H(\ell_1) = C_n \times C_m$ ($n,m \in \mathbb{N}$), 
then 
$$
\mathcal{K}(H, \varphi, L) = \mathcal{K}(\mathbb{T}^2 , \varphi, \mathbb{Z}_q(2)) = \mathbb{T}^2 \cup \mathbb{T}^2
$$
is a commutative hypergroup. 

Obviously $\mathcal{K}(\mathbb{T}^2, \varphi, \mathbb{Z}_q(2))$ is Pontryagin 
and 
$$
\hat{\mathcal{K}}(\mathbb{T}^2, \varphi, \mathbb{Z}_q(2)) = 
\mathcal{K}(\mathbb{Z}_q(2), \hat{\varphi}, \mathbb{Z}^2), 
$$
where the dual field $\hat{\varphi}$ of $\varphi$ is given as 
$$
\hat{\varphi} :  \mathbb{Z} \times \mathbb{Z}  \ni k \longmapsto \hat{\varphi}(k) \subset  \mathbb{Z}_q(2)
$$
with 
\begin{eqnarray*}
\hat{\varphi}(k)=
\left\{ 
\begin{array}{ll}
\{\ell_0\} &  ~~~{\rm for}~~k \in n \mathbb{Z} \times m \mathbb{Z}  \\
\mathbb{Z}_q(2) &  ~~~{\rm otherwise}. \\
\end{array} 
\right.
\end{eqnarray*}

\bigskip

{\bf Example 4.3 } If  
$H = \mathbb{T} $,  
$L = \mathbb{Z}_q(3) = \{\ell_0, \ell_1,\ell_2 \}$ (see [KTY]) and 
$\varphi : \mathbb{Z}_q(3) \ni \ell \longmapsto H(\ell) \subset H$ 
with $\varphi(\ell_0) = H(\ell_0) = \{ 1 \}$, 
$\varphi(\ell_1) = H(\ell_1) = C_n $ and $\varphi(\ell_2) = H(\ell_2) = C_n $, 
then 
$$
\mathcal{K}(H, \varphi, L) = \mathcal{K}(\mathbb{T} , \varphi, \mathbb{Z}_q(3)) 
= \mathbb{T} \cup \mathbb{T} \cup \mathbb{T}
$$
is a commutative hypergroup. 

Obviously $\mathcal{K}(\mathbb{T}, \varphi, \mathbb{Z}_q(3))$ is Pontryagin 
and 
$$
\hat{\mathcal{K}}(\mathbb{T}, \varphi, \mathbb{Z}_q(3)) = 
\mathcal{K}(\mathbb{Z}_q(3), \hat{\varphi}, \mathbb{Z}), 
$$
where the dual field $\hat{\varphi}$ of $\varphi$ is given as 
$$
\hat{\varphi} :  \mathbb{Z} \ni k \longmapsto \hat{\varphi}(k) \subset \mathbb{Z}_q(3)  
$$
with 
\begin{eqnarray*}
\hat{\varphi}(k)=
\left\{ 
\begin{array}{ll}
\{\ell_0\} &  ~~~{\rm for}~~k \in n \mathbb{Z}   \\
\mathbb{Z}_q(3) &  ~~~{\rm otherwise}. \\
\end{array} 
\right.
\end{eqnarray*}

\bigskip

{\bf Example 4.4 } If  
$H = \mathbb{T} \times \mathbb{T} = \mathbb{T}^2$,  
$L = \mathbb{Z}_q(2) = \{\ell_0, \ell_1\}$ and 
$\varphi : \mathbb{Z}_q(2) \ni \ell \longmapsto H(\ell) \subset H$ 
with $\varphi(\ell_0) = H(\ell_0) = \{ (1,1) \}$, 
$\varphi(\ell_1) = H(\ell_1) = C_n \times \mathbb{T}$ ($n \in \mathbb{N}$), 
then 
$$
\mathcal{K}(H, \varphi, L) = \mathcal{K}(\mathbb{T}^2 , \varphi, \mathbb{Z}_q(2)) = \mathbb{T}^2 \cup \mathbb{T}
$$
is a commutative hypergroup. 

Obviously $\mathcal{K}(\mathbb{T}^2, \varphi, \mathbb{Z}_q(2))$ is Pontryagin 
and 
$$
\hat{\mathcal{K}}(\mathbb{T}^2, \varphi, \mathbb{Z}_q(2)) = 
\mathcal{K}(\mathbb{Z}_q(2), \hat{\varphi}, \mathbb{Z}^2), 
$$
where the dual field $\hat{\varphi}$ of $\varphi$ is given as 
$$
\hat{\varphi} :  \mathbb{Z} \times  \mathbb{Z} 
\ni k \longmapsto \hat{\varphi}(k) \subset \mathbb{Z}_q(2)
$$
with 
\begin{eqnarray*}
\hat{\varphi}(k)=
\left\{ 
\begin{array}{ll}
\{\ell_0 \} &  ~~~{\rm for}~~k \in n \mathbb{Z} \times \{ 0 \}  \\
\mathbb{Z}_q(2) &  ~~~{\rm otherwise}. \\
\end{array} 
\right.
\end{eqnarray*}

\bigskip

{\bf Example 4.5 } Let   
$H = \mathcal{K}^{\alpha}(\mathbb{T})  = [-1, 1]  $,  
$L = \mathbb{Z}_q(2) = \{\ell_0, \ell_1\}$ and 
$\varphi : \mathbb{Z}_q(2) \ni \ell \longmapsto H(\ell) \subset H$ 
with $\varphi(\ell_0) = H(\ell_0) = \{ 1 \}$, 
$\varphi(\ell_1) = H(\ell_1) =  \mathcal{K}^{\alpha}(C_n) $ ($n \in \mathbb{N}$), 
 where $\mathcal{K}^{\alpha}(\mathbb{T})  = [-1, 1]$ is the orbital hypergroup defined by the 
action of  $\mathbb{Z}_2 = \{ e, g \} \ ( g^2=e) $
such that $\alpha_g(z) = \bar{z}$ on $\mathbb{T}$.
Then 
$$
\mathcal{K}(H, \varphi, L) = \mathcal{K}(\mathcal{K}^{\alpha}(\mathbb{T}), \varphi, \mathbb{Z}_q(2)) = [-1, 1] \cup [-1, 1] 
$$
is a commutative hypergroup. 

Obviously $\mathcal{K}(\mathcal{K}^{\alpha}(\mathbb{T}),  \varphi, \mathbb{Z}_q(2))$ is Pontryagin 
and 
$$
\hat{\mathcal{K}}(\mathcal{K}^{\alpha}(\mathbb{T}), \varphi, \mathbb{Z}_q(2)) = 
\mathcal{K}(\mathbb{Z}_q(2), \hat{\varphi}, \mathcal{K}^{\alpha}(\mathbb{Z})), 
$$
where $\mathcal{K}^{\alpha}(\mathbb{Z})  = \{ 0,1,2, \cdots, n, \cdots \}$ is also 
the orbital hypergroup by the 
action of  $\mathbb{Z}_2 = \{ e, g \} \ ( g^2=e) $
such that $\alpha_g(n) = -n $ on $\mathbb{Z}$ and the dual field $\hat{\varphi}$ of $\varphi$ is given as 
$$
\hat{\varphi} : \mathcal{K}^{\alpha}(\mathbb{Z}) \ni k \longmapsto \hat{\varphi}(k) \subset 
\mathbb{Z}_q(2)  
$$
with 
\begin{eqnarray*}
\hat{\varphi}(k)=
\left\{ 
\begin{array}{ll}
\{\ell_0\} &  ~~~{\rm for}~~k \in \mathcal{K}^{\alpha}(n \mathbb{Z})  \\
\mathbb{Z}_q(2) &  ~~~{\rm otherwise}. \\
\end{array} 
\right.
\end{eqnarray*}

\bigskip

{\bf Example 4.6 } Let $A$ be a commutative strong hypergroup and $C$ a compact strong hypergroup. If  
$H = A \times C$,  
$L = \mathbb{Z}_q(2) = \{\ell_0, \ell_1\}$ and 
$\varphi : \mathbb{Z}_q(2) \ni \ell \longmapsto H(\ell) \subset H$ 
with $\varphi(\ell_0) = H(\ell_0) = \{ h_0 \}$, 
$\varphi(\ell_1) = H(\ell_1) = C$,  
then 
$$
\mathcal{K}(H, \varphi, L) = \mathcal{K}(A \times C, \varphi, \mathbb{Z}_q(2)) = (A \times C) \cup A
$$
is the commutative hypergroup $A \times (C \vee \mathbb{Z}_q(2))$ of strong type  
and 
$$
\widehat{\mathcal{K}}(A \times C, \varphi, \mathbb{Z}_q(2)) = \mathcal{K}( \mathbb{Z}_q(2), \hat{\varphi}, \hat{A} \times \hat{C}) = \hat{A} \times ( \mathbb{Z}_q(2) \vee \hat{C}).
$$
In fact, the dual field $\hat{\varphi}$ of $\varphi$ is given as 
$$
\hat{\varphi} : \hat{A} \times \hat{C} \ni (\chi, \rho) \longmapsto \hat{\varphi} (\chi, \rho) \subset 
\mathbb{Z}_q(2)  
$$
with 
\begin{eqnarray*}
\hat{\varphi} (\chi, \rho)  =
\left\{ 
\begin{array}{ll}
\{\ell_0\} &  ~~~{\rm for}~~ (\chi, \rho_0) \in \hat{A} \times \hat{C}  \\
\mathbb{Z}_q(2) &  ~~~{\rm otherwise}, \\
\end{array} 
\right.
\end{eqnarray*}
where $\rho_0$ is unit of $\hat{C}$.

\bigskip

{\bf Example 4.7 } Let $A$ be a commutative hypergroup and $C$ a compact hypergroup. If  
$H = A \times C$,  
$L = \mathbb{Z}_q(3) = \{\ell_0, \ell_1, \ell_2\}$ and 
$\varphi : \mathbb{Z}_q(3) \ni \ell \longmapsto H(\ell) \subset H$ 
with $\varphi(\ell_0) = H(\ell_0) = \{ h_0 \}$, 
$\varphi(\ell_1) = H(\ell_1) = C_0$, $ \varphi(\ell_2) = H(\ell_2)  = C_0$, where 
$C_0$ is a closed subhypergroup of $C$,  
then 
$$
\mathcal{K}(H, \varphi, L) = \mathcal{K}(A \times C, \varphi, \mathbb{Z}_q(3))
=(A \times C) \cup (A \times Q) \cup (A \times Q)
$$
is the commutative hypergroup $A \times (S(Q \times \mathbb{Z}_q(3) : Q \rightarrow C))$,  
where $Q = C/C_0$.

\bigskip

\textbf{\large 5. Applications of the theorems }

\medskip

We assume $H$ to be a (not necessarily compact) commutative hypergroup of 
strong type such that $\hat{H}$ is a commutative hypergroup with unit character $\chi_0$. 
Let $H_0$ be a closed subhypergroup of strong type of $H$, where the annihilator 
$H_0^\bot$ in $\hat{H}$ is a compact subhypergroup of $\hat{H}$. 
For $\tau \in \widehat{H_0}$ we consider the set 
$$
A(\tau) = \{\chi \in \hat{H} : {\rm res}_{H_0}^H\ \chi = \tau \}.
$$
As usual $\omega_{H_0^\bot}$ denotes the normalized Haar measure of $H_0^\bot$. 
Then there exists a unique $H_0^\bot$-invariant probability measure $\mu_{A(\tau)}$ 
which is given by 
$$
\mu_{A(\tau)} = ch(\tilde{\tau}) \cdot \omega_{H_0^\bot}.
$$
for some $\tilde{\tau} \in A(\tau)$. 

\medskip
\noindent
{\bf Definition } (see [HKY]) 

\medskip
(i) For $\tau \in \widehat{H_0}$ the character of $\tau$ induced from 
$H_0$ to $H$ is defined by 
$$
{\rm ind}_{H_0}^H ch(\tau) := \mu_{A(\tau)}. 
$$

(ii) For $\tau_i, \tau_j \in \widehat{H_0}$, $ch(\tau_i) \cdot ch(\tau_j)$ 
is decomposed on $\widehat{H_0}$ in the form 
$$
ch(\tau_i) \cdot  ch(\tau_j) = \int_C ch(\tau) \nu (d \tau), 
$$ 
where $\nu$ is a probability measure on $\widehat{H_0}$ and 
$$
C := {\rm supp}(\nu) = {\rm supp}(ch(\tau_i) \cdot ch(\tau_j))
$$
is compact. 

\medskip

Then we introduce 
$$
{\rm ind}_{H_0}^H (ch(\tau_i) \cdot ch(\tau_j)) 
:= \int_C {\rm ind}_{H_0}^H ch(\tau) \nu(d \tau). 
$$
The subsequent simple facts play an essential role in the upcoming discussion.

\medskip
\noindent
{\bf Lemma 5.1 } (see [HKY])

\medskip
(i) For $\tau \in \widehat{H_0}$ 
$$
{\rm res}_{H_0}^H ({\rm ind}_{H_0}^H ch(\tau)) = ch(\tau). 
$$

\medskip
(ii) For $\pi_i, \pi_j \in \hat{H}$
$$
{\rm res}_{H_0}^H (ch(\pi_i) \cdot ch(\pi_j)) = 
({\rm res}_{H_0}^H ch(\pi_i)) \cdot ({\rm res}_{H_0}^H ch(\pi_j)). 
$$

\medskip
(iii) For $\pi \in \hat{H}$ and $\tau_i, \tau_j \in \widehat{H_0}$
$$
{\rm ind}_{H_0}^H ( ({\rm res}_{H_0}^H ch(\pi)) \cdot ch(\tau_i) \cdot ch(\tau_j))
 = ch(\pi) \cdot {\rm ind}_{H_0}^H (ch(\tau_i) \cdot ch(\tau_j) )
$$

\medskip

(iv) For $\tau_i, \tau_j \in \widehat{H_0}$
$$
{\rm res}_{H_0}^H ({\rm ind}_{H_0}^H (ch(\tau_i) \cdot  ch(\tau_j))) = 
ch(\tau_i) ch(\tau_j). 
$$

\medskip
\noindent
{\bf Proof } (i) and (ii) are clear. 

\medskip
(iii) It is easy to check that 
$$
{\rm ind}_{H_0}^H (ch(\tau_i) \cdot ch(\tau_j)) 
= ch(\tilde{\tau_i}) \cdot ch(\tilde{\tau_j}) \cdot \omega_{H_0^\bot}
$$
for $\tilde{\tau_i} \in A(\tau_i)$ and $\tilde{\tau_j} \in A(\tau_j)$. Then 
we see that 
$$
ch(\pi) \cdot {\rm ind}_{H_0}^H (ch(\tau_i) \cdot ch(\tau_j)) 
= ch(\pi) \cdot  ch(\tilde{\tau_i}) \cdot ch(\tilde{\tau_j}) \cdot \omega_{H_0^\bot}
$$
and 
\begin{align*}
{\rm ind}_{H_0}^H ( ({\rm res}_{H_0}^H ch(\pi)) \cdot ch(\tau_i) \cdot ch(\tau_j))  
= ch(\pi) \cdot  ch(\tilde{\tau_i}) \cdot ch(\tilde{\tau_j}) \cdot \omega_{H_0^\bot}. 
\end{align*}

\medskip
(iv) For $\tau_i, \tau_j \in \widehat{H_0}$
\begin{align*}
{\rm res}_{H_0}^H ({\rm ind}_{H_0}^H (ch(\tau_i) \cdot  ch(\tau_j))) 
&= {\rm res}_{H_0}^H \left(\int_C  {\rm ind}_{H_0}^H ch(\tau) \nu (d\tau) \right)\\
&= \int_C {\rm res}_{H_0}^H ({\rm ind}_{H_0}^H ch(\tau)) \nu (d\tau)\\
&= \int_C ch(\tau) \nu (d \tau)\\
&= ch(\tau_i) \cdot ch(\tau_j).
\end{align*}

\bigskip

\noindent
{\bf Remark}  \ A pair $(H,H_0)$ for a commutative hypergroup of strong type 
is always an admissible hypergroup pair in the sense of [HKTY2] by Lemma 5.1.

\bigskip

\noindent
{\bf Definition} (see [HKTY2])\ On the space  $$
\mathcal{K}(\hat{H} \cup \widehat{H_0},\mathbb{Z}_q(2))
:= \{(ch(\pi), \circ), (ch(\tau), \bullet) : \pi \in \hat{H}, \tau \in \widehat{H_0}\}  
$$ 
we define the convolution $*$ by the following properties:

\begin{enumerate}
\item
$(ch(\pi_i),\circ)*(ch(\pi_j),\circ):=
(ch(\pi_i) \cdot ch(\pi_j),\circ),$

\medskip

\item
$(ch(\pi),\circ)*(ch(\tau),\bullet):=
(({\rm res}_{H_0}^H ch(\pi)) \cdot ch(\tau),\bullet),$

\medskip

\item
$(ch(\tau),\bullet)*(ch(\pi),\circ):=
(ch(\tau) \cdot ({\rm res}_{H_0}^H ch(\pi)),\bullet),$

\medskip

\item
$
(ch(\tau_i),\bullet)*(ch(\tau_j),\bullet):=
q({\rm ind}_{H_0}^H(ch(\tau_i) \cdot ch(\tau_j)),\circ) +(1-q)(ch(\tau_i) \cdot ch(\tau_j),\bullet).
$
\end{enumerate}

\medskip

\medskip
\noindent
{\bf Definition } Given $\mathbb{Z}_q(2) = \{\ell_0, \ell_1\}$ $(0 < q \leq 1)$ 
we introduce the set $\mathcal{K}(\hat{H}, \varphi, \mathbb{Z}_q(2))$ via 
the hyperfield  $\varphi : \mathbb{Z}_q(2) \ni \ell \mapsto \hat{\varphi}(\ell) \subset \hat{H}$ given by 
\begin{eqnarray*}
\varphi(\ell)=
\left\{ 
\begin{array}{ll}
\{\chi_0\} &  ~{\rm if}~~\ell = \ell_0  \\
H_0^\bot &  ~{\rm if}~~\ell = \ell_1 \\
\end{array} 
\right.
\end{eqnarray*}
as in section 3. 

\medskip
\noindent
Then we have 

\bigskip

\medskip
\noindent
{\bf Theorem 5.2 } Let $H$ be a commutative hypergroup of strong type and 
$H_0$ a closed subhypergroup of $H_0$
such that $H_0^{\bot}$ is compact in $\hat{H}$. Then  
$\mathcal{K}(\hat{H} \cup \widehat{H_0}, \mathbb{Z}_q(2))$ is a commutative 
hypergroup and 
$$
\mathcal{K}(\hat{H} \cup \widehat{H_0}, \mathbb{Z}_q(2)) \cong 
\mathcal{K}(\hat{H}, \varphi, \mathbb{Z}_q(2)). 
$$

\medskip
\noindent
{\bf Proof }  
In order to show that $\mathcal{K}(\hat{H} \cup \widehat{H_0}, \mathbb{Z}_q(2))$ 
is a hypergroup we should check the following associativity relations. 
For $\pi_i, \pi_j, \pi_k, \pi \in \hat{H}$ and $\tau_i, \tau_j, \tau_k, \tau \in \widehat{H_0}$ 

\begin{align*}
&(A1)~~((ch(\pi_i), \circ) * (ch(\pi_j), \circ)) * (ch(\pi_k), \circ) 
=  (ch(\pi_i), \circ) * ((ch(\pi_j), \circ) * (ch(\pi_k), \circ)).\\
&(A2)~~((ch(\pi_i), \circ) * (ch(\pi_j), \circ)) * (ch(\tau), \bullet) = 
(ch(\pi_i), \circ) * ((ch(\pi_j), \circ) * (ch(\tau), \bullet)).\\
&(A3)~~((ch(\pi), \circ) * (ch(\tau_i), \bullet)) * (ch(\tau_j), \bullet) = 
(ch(\pi), \circ) * ((ch(\tau_i), \bullet) * (ch(\tau_j), \bullet)).\\
&(A4)~~((ch(\tau_i), \bullet) * (ch(\tau_j), \bullet)) * (ch(\tau_k), \bullet) = 
(ch(\tau_i), \bullet) * ((ch(\tau_j), \bullet) * (ch(\tau_k), \bullet)).
\end{align*}
However these relations are shown in a similar way to the proof of Proposition 3.6 
in our paper [HKTY2] combined with the above Lemma 5.1 so that we omit the details.   
It is easy to check the remaining axioms of a hypergroup 
for $\mathcal{K}(\hat{H} \cup \widehat{H_0}, \mathbb{Z}_q(2))$.  
The desired conclusion is obtained.

\medskip

Next we introduce an isomorphism 
$\psi : \mathcal{K}(\hat{H} \cup \widehat{H_0}, \mathbb{Z}_q(2))
\longrightarrow \mathcal{K}(\hat{H}, \varphi, \mathbb{Z}_q(2))$ by 
\begin{align*}
\psi((ch(\pi), \circ)) = ch(\pi) \otimes \varepsilon_{\ell_0},~~~
\psi((ch(\tau), \bullet)) = (ch(\tilde{\tau}) \cdot \omega_{H_0^\bot}) \otimes \varepsilon_{\ell_1}. 
\end{align*}
It is easy to see that $\psi$ is bijective. We only show that $\psi$ is  
homomorphic. 
\begin{align*}
1.~~\psi((ch(\pi_i),\circ)*(ch(\pi_j),\circ)) 
&= \psi((ch(\pi_i) \cdot ch(\pi_j),\circ)) \\
&= (ch(\pi_i) \cdot ch(\pi_j)) \otimes \varepsilon_{\ell_0}\\
&= (ch(\pi_i) \otimes \varepsilon_{\ell_0}) \circ (ch(\pi_j) \otimes \varepsilon_{\ell_0})\\
&= \psi((ch(\pi_i),\circ)) \circ \psi((ch(\pi_j),\circ)), 
\end{align*}
\begin{align*}
2.~~ \psi((ch(\pi),\circ)*(ch(\tau),\bullet)) 
&= \psi((({\rm res}_{H_0}^H ch(\pi)) \cdot ch(\tau),\bullet))\\
&= (ch(\pi) \cdot ch(\tilde{\tau}) \cdot \omega_{H_0}^\bot) \otimes \varepsilon_{\ell_1}\\
&= (ch(\pi) \otimes \varepsilon_{\ell_0}) \circ 
(ch(\tilde{\tau}) \cdot \omega_{H_0}^\bot \otimes \varepsilon_{\ell_1})\\
&= \psi((ch(\pi),\circ)) \circ \psi((ch(\tau),\bullet))
\end{align*}

$$
\hspace{-10mm}3.~~\psi((ch(\tau),\bullet)*(ch(\pi),\circ)) =
\psi((ch(\tau),\bullet)) \circ \psi((ch(\pi),\circ))
$$ 
is obtained similarly. 
\begin{align*}
\hspace{15mm}4.~~ 
&\psi((ch(\tau_i),\bullet)*(ch(\tau_j),\bullet))\\
= & \psi( q({\rm ind}_{H_0}^H(ch(\tau_i) \cdot ch(\tau_j)),\circ)
+(1-q)(ch(\tau_i) \cdot ch(\tau_j),\bullet))\\
= & q\psi(({\rm ind}_{H_0}^H(ch(\tau_i) \cdot ch(\tau_j)),\circ)) 
+(1-q)\psi((ch(\tau_i) \cdot ch(\tau_j),\bullet))\\
= & q((ch(\tilde{\tau_i}) \cdot ch(\tilde{\tau_j}) \cdot \omega_{H_0^\bot}) \otimes \varepsilon_{\ell_0}) 
+(1-q)((ch(\tilde{\tau_i}) \cdot ch(\tilde{\tau_j}) \cdot \omega_{H_0^\bot}) \otimes \varepsilon_{\ell_1})\\
=&(ch(\tilde{\tau_i}) \cdot ch(\tilde{\tau_j}) \cdot \omega_{H_0^\bot}) \otimes 
(q \varepsilon_{\ell_0} + (1-q)\varepsilon_{\ell_1}) \\
=&(ch(\tilde{\tau_i}) \cdot ch(\tilde{\tau_j}) \cdot \omega_{H_0^\bot}) \otimes 
(\varepsilon_{\ell_1} \bullet \varepsilon_{\ell_1}) \\
=& ((ch(\tilde{\tau_i}) \cdot \omega_{H_0^\bot}) \otimes \varepsilon_{\ell_1}) \circ 
((ch(\tilde{\tau_j}) \cdot \omega_{H_0^\bot}) \otimes \varepsilon_{\ell_1})\\
=& \psi((ch(\tau_i), \bullet)) \circ \psi((ch(\tau_j), \bullet)). 
&\hspace{-15mm}\text{[Q.E.D.]}
\end{align*}

\medskip
Now let $H$ be a discrete commutative hypergroup of Pontryagin type and 
$H_0$ a closed subhypergroup of $H$. Since $\hat{H}$ is a compact hypergroup, 
$H_0^\bot$ is a compact subhypergroup. 
Then $\mathcal{K}(\hat{H}, \varphi, \mathbb{Z}_q(2))$ is 
defined via the  hyperfield $\varphi$  given by 
\begin{eqnarray*}
\varphi(\ell)=
\left\{ 
\begin{array}{ll}
\{\chi_0\} &  ~{\rm if}~~\ell = \ell_0  \\
H_0^\bot &  ~{\rm if}~~\ell = \ell_1. \\
\end{array} 
\right.
\end{eqnarray*}
The dual field $\hat{\varphi}$ of $\varphi$ is the field 
$$
\hat{\varphi} : H \ni h \longmapsto \hat{\varphi}(h) \subset \mathbb{Z}_q(2)
$$
with 
\begin{eqnarray*}
\hat{\varphi}(h)=
\left\{ 
\begin{array}{ll}
\{\ell_0\} &  ~{\rm if}~~h \in H_0  \\
\mathbb{Z}_q(2) &  ~{\rm otherwise}.\\
\end{array} 
\right.
\end{eqnarray*}
Applying Theorem 3.5 together with Theorem 5.2 we obtain

\bigskip
\noindent
{\bf Theorem 5.3 } Let $H$ be a discrete commutative hypergroup of Pontryagin type and 
$H_0$ a closed subhypergroup of $H$ such that $H_0^{\bot}$ is compact in $\hat{H}$. Then
$$
\hat{\mathcal{K}}(\hat{H} \cup \widehat{H_0}, \mathbb{Z}_q(2)) \cong  
\mathcal{K}(\mathbb{Z}_q(2), \hat{\varphi}, H). 
$$

\bigskip

\bigskip

\textbf{\large 6. Examples of hypergroup duals of $\mathcal{K}(\hat{H} \cup \widehat{H_0}, \mathbb{Z}_q(2))$}

\medskip

\medskip

\noindent
{\bf Example 6.1 } If 
$H = \mathbb{Z} $ and $H_0 = n\mathbb{Z} $ \  ($n \in \mathbb{N}$), then 
$$\mathcal{K}(\hat{\mathbb{Z}} \cup \widehat{n\mathbb{Z}}, \mathbb{Z}_q(2)) = 
\mathcal{K}(\mathbb{T}, \varphi, \mathbb{Z}_q(2)) = \mathbb{T} \cup \mathbb{T} \ \ 
{\rm (in \ Example \ 4.1)}
$$
and
$$
\widehat{\mathcal{K}}(\hat{\mathbb{Z}} \cup \widehat{n\mathbb{Z}}, \mathbb{Z}_q(2)) = 
\mathcal{K}(\mathbb{Z}_q(2), \hat{\varphi}, \mathbb{Z}) \ \ 
{\rm (in \ Example \ 4.1)}. 
$$

\bigskip

\noindent
{\bf Example 6.2 } If 
$H = \mathbb{Z}^2 =  \mathbb{Z} \times \mathbb{Z}$ and $H_0 = n\mathbb{Z} \times m\mathbb{Z} $ 
\  ($n, m \in \mathbb{N}$), then 
$$\mathcal{K}(\widehat{\mathbb{Z}^2} \cup (\widehat{n\mathbb{Z} \times m\mathbb{Z}}), \mathbb{Z}_q(2)) = 
\mathcal{K}(\mathbb{T}^2 , \varphi, \mathbb{Z}_q(2)) = \mathbb{T}^2 \cup \mathbb{T}^2
\ \ 
{\rm (in \ Example \ 4.2)} 
$$
and
$$\widehat{\mathcal{K}}(\widehat{\mathbb{Z}^2} \cup (\widehat{n\mathbb{Z} \times m\mathbb{Z}}), \mathbb{Z}_q(2)) = 
\mathcal{K}(\mathbb{Z}_q(2),  \hat{\varphi}, \mathbb{Z}^2) 
\ \ 
{\rm (in \ Example \ 4.2)}. 
$$

\bigskip

\bigskip

\noindent
{\bf Example 6.3 } If 
$H = \mathbb{Z}^2 =  \mathbb{Z} \times \mathbb{Z}$ and $H_0 = n\mathbb{Z} \times \{0\} \cong n\mathbb{Z} $ 
\  ($n \in \mathbb{N}$), then 
$$\mathcal{K}(\widehat{\mathbb{Z}^2} \cup (\widehat{n\mathbb{Z} }), \mathbb{Z}_q(2)) = 
\mathcal{K}(\mathbb{T}^2 , \varphi, \mathbb{Z}_q(2)) = \mathbb{T}^2 \cup \mathbb{T}
\ \ 
{\rm (in \ Example \ 4.4)} 
$$
and
$$\widehat{\mathcal{K}}(\widehat{\mathbb{Z}^2} \cup (\widehat{n\mathbb{Z} }), \mathbb{Z}_q(2)) = 
\mathcal{K}(\mathbb{Z}_q(2), \hat{\varphi}, \mathbb{Z}^2) 
\ \ 
{\rm (in \ Example \ 4.4)}. 
$$

\bigskip

\noindent
{\bf Example 6.4 } If 
$H = \mathcal{K}^{\alpha}(\mathbb{Z})$ and $H_0 =  \mathcal{K}^{\alpha}(n\mathbb{Z}) $ 
\  ($n \in \mathbb{N}$), then 
$$\mathcal{K}(\widehat{ \mathcal{K}^{\alpha}(\mathbb{Z})} \cup (\widehat{ \mathcal{K}^{\alpha}(n\mathbb{Z})}), \mathbb{Z}_q(2)) = 
\mathcal{K}( \mathcal{K}^{\alpha}(\mathbb{T}), \varphi, \mathbb{Z}_q(2)) 
\ \ 
{\rm (in \ Example \ 4.5)} 
$$
and
$$\widehat{\mathcal{K}}(\widehat{ \mathcal{K}^{\alpha}(\mathbb{Z})} \cup 
(\widehat{ \mathcal{K}^{\alpha}(n\mathbb{Z})}), \mathbb{Z}_q(2))  = 
\mathcal{K}( \mathbb{Z}_q(2), \hat{\varphi}, \mathcal{K}^{\alpha}(\mathbb{Z})) 
\ \ 
{\rm (in \ Example \ 4.5)}. 
$$

\bigskip

\bigskip

\noindent
{\bf Example 6.5 } Let $B$ be a commutative Pontryagin hypergroup and $D$ a discrete 
commutative Pontryagin hypergroup. If  
$H = B \times D$ and   
$H_0 = B$, then 
$$\mathcal{K}(\widehat{B \times D} \cup \widehat{B}, 
\mathbb{Z}_q(2)) = 
\mathcal{K}( \hat{B} \times \hat{D}, \varphi, \mathbb{Z}_q(2))  = 
\hat{B} \times (\hat{D} \vee \mathbb{Z}_q(2))
\ \ 
{\rm (in \ Example \ 4.6)} 
$$
and 
$$
\widehat{\mathcal{K}}(\widehat{B \times D} \cup \widehat{B}, \mathbb{Z}_q(2)) = \mathcal{K}( \mathbb{Z}_q(2), \hat{\varphi}, 
B \times D) = B \times ( \mathbb{Z}_q(2) \vee D))
\ \ 
{\rm (in \ Example \ 4.6)}. 
$$

\bigskip

\textbf{Addresses }

\medskip

Herbert Heyer : Universit\"{a}t T\"{u}bingen

Mathematisches Institut

Auf der Morgenstelle 10

72076, T\"{u}bingen

Germany

\medskip

e-mail : herbert.heyer@uni-tuebingen.de

\bigskip

Satoshi Kawakami : Nara University of Education

Department of Mathematics

Takabatake-cho

Nara, 630-8528

Japan

\medskip

e-mail : kawakami@nara-edu.ac.jp

\bigskip
Tatsuya Tsurii : Osaka Prefecture University 

Graduate School of Science

1-1 Gakuen-cho, Nakaku, Sakai 

Osaka, 599-8531

Japan

\medskip

e-mail : dw301003@edu.osakafu-u.ac.jp

\bigskip

Satoe Yamanaka : Nara Women's University 

Faculty of Science 

Kita-uoya-higashi-machi, 

Nara, 630-8506

Japan

\medskip

e-mail : s.yamanaka516@gmail.com

\end {document}